\documentclass[10pt]{amsart}
\input{epsf}
\usepackage{amssymb,latexsym}
\topskip  \usepackage{euscript,color}
\usepackage{graphpap}
\usepackage{tikz}

\newtheorem{theorem}{Theorem}

\newtheorem{corollary}[theorem]{Corollary}

\newtheorem{lemma}[theorem]{Lemma}

\begin{document}

\title{Modelling x-ray tomography using integer compositions}
\maketitle

\begin{center}Aubrey Blecher\\
School of Mathematics, University of Witwatersrand, Johannesburg, SA

{\tt Aubrey.Blecher@wits.ac.za}
\end{center}\begin{center}

Toufik Mansour\\
Department of Mathematics, University of Haifa, 31905 Haifa, Israel

{\tt tmansour@univ.haifa.ac.il}
\end{center}

\begin{abstract}
The x-ray process is modelled using integer compositions as a two dimensional analogue of the object being x-rayed, where the examining rays are modelled by diagonal lines with equation $x-y=n$ for non negative integers $n$. This process is essentially parameterised by the degree to which the x-rays are contained inside a particular composition. So, characterising the process translates naturally to obtaining a generating function which tracks the number of "staircases" which are contained inside arbitrary integer compositions of $n$. More precisely, we obtain a generating function which counts the number of times the staircase $1^+2^+3^+\cdots m^+$ fits inside a particular composition. The main theorem establishes this generating function
\begin{equation*}
F= \dfrac {k_{m}-\frac {qx^{m}y}{1-x}k_{m-1}}{(1-q)x^{\binom {m+1}{2}}\left(\frac{y}{1-x}\right)^{m}+\frac{1-x-xy}{1-x}\left(k_{m}-\frac{qx^{m}y}{1-x}k_{m-1}\right)}.
\end{equation*}
where
\begin{equation*}
k_{m}=\sum_{\j=0}^{m-1}x^{mj-\binom {j}{2}}\left(\frac {y}{1-x}\right)^{j}.
\end{equation*}
Here $x$ and $y$ respectively track the composition size and number of parts, whilst $q$ tracks the number of such staircases contained.
\end{abstract}



\section{Introduction}
In several recent papers the notion of integer compositions of $n$ (represented as the associated bargraph) have been used to model certain problems in physics. See for example \cite{LBG,FATP,ESAW,ESAP} where bargraphs are a representation of a polymer at an adsorbing wall subject to several forces.

In a paper by a current author et al (see \cite{BMPS}), the x-ray process was modelled using permutation matrices as a two dimensional analogue of the object being x-rayed, where the examining rays are modelled by diagonal lines with equation $x+y=n$ for positive integers $n$. The current paper is based instead on integer compositions as the object analogue and where the examining rays are represented by equation $x-y=n$ for non negative integers $n$. Since this model is essentially parameterized by the degree to which the x-rays are contained inside an arbitrary composition, it translates naturally to obtaining a generating function which tracks the number of "staircases" which are contained inside particular integer compositions of $n$. More precisely, we will obtain a generating function which counts (with the exponent $s$ of $q$ as tracker) the number of times the staircase $1^+2^+3^+\cdots m^+$ ($m$ fixed) fits inside particular compositions. So the term of our generating function $n(a,b,s)x^ay^bq^s$ indicates that there are in total $n(a,b,s)$ compositions of $a$ with $b$ parts in which the staircases $1^+2^+3^+\cdots m^+$ occurs exactly $s$ times.

\subsection{Definitions} A
{\it composition of a positive integer} $n$ is a sequence of $k$ positive integers $a_{1},a_{2},\cdots a_{k}$, each called a part such that $n=\sum_ {\i=1}^{k}a_{i}$; A {\it staircase $1^+2^+3^+\cdots m^+$} is a word with m sequential parts from left to right where for $1 \le i \le m$ the $i$th part $\ge i$.\\\\
See for example the staircase in Figure \ref{staircase} below.
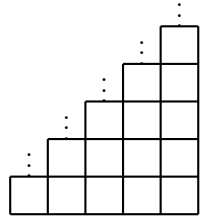
\begin{figure}[h!]
\begin{center}
\begin{tikzpicture}[scale=0.5]
\draw [thick] (0,0)--(5,0)(0,1)--(5,1)(1,2)--(5,2)(2,3)--(5,3)(3,4)--(5,4)(4,5)--(5,5);
\draw [thick] (0,0)--(0,1)(1,0)--(1,2)(2,0)--(2,3)(3,0)--(3,4)(4,0)--(4,5)(5,0)--(5,5);
\draw (0.5,1.5)node[]{\vdots};\draw (1.5,2.5)node[]{\vdots};\draw (2.5,3.5)node[]{\vdots};\draw (3.5,4.5)node[]{\vdots};\draw (4.5,5.5)node[]{\vdots};
\end{tikzpicture}\caption{The staircase $1^+2^+3^+4^+5^+$}\label{staircase}
\end{center}
\end{figure}

Much recent work has been done on various statistics relating to compositions. See, for example, \cite{FASd,SPDC,SPC} and \cite{Mb1} and references therein.

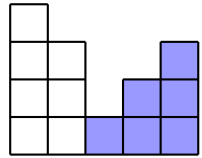
\begin{figure}[h!]
\begin{center}
\begin{tikzpicture}[scale=0.5]
\fill[blue!40!white] (2,0)--(2,1)--(3,1)--(3,2)--(4,2)--(4,3)--(5,3)--(5,0);
\draw [thick] (0,0)--(5,0)(0,1)--(5,1)(0,2)--(2,2)(3,2)--(5,2)(0,3)--(2,3)(4,3)--(5,3)(0,4)--(1,4);
\draw [thick] (0,0)--(0,4)(1,0)--(1,4)(2,0)--(2,3)(3,0)--(3,2)(4,0)--(4,3)(5,0)--(5,3);
\end{tikzpicture}\caption{The composition $4+3+1+2+3$ containing one staircase $1^+2^+3^+$ (coloured) and three $1^+2^+$ staircases}\label{composition}
\end{center}
\end{figure}

A particular composition may be represented as a bargraph (see \cite{Mb1} and \cite{LBG}). For example the composition $4+3+1+2+3$ of $13$ represented in Figure \ref{composition} as a bargraph, contains exactly one $1^+2^+3^+$ staircase, three $1^+2^+$ staircases and five $1^+$ staircases. It contains no others.

In this paper, compositions (ie their associated bargraphs) are the analogue for a (2-dimensional) object to be x-rayed (as explained above). Across all possible compositions, the shapes are parameterized in a generating function by a marker variable $q$ which tracks the number of $1^+2^+3^+\cdots m^+$ staircases (again with $m$ fixed) that fit inside a composition. The generating function in question is defined as
\begin{equation}\label{f36}
F=\sum_{n\geq1; b\geq1; s\geq 0}n(a,b,c)x^{a}y^{b}q^{s},
\end{equation}
where $n(a,b,s)$ is the number of compositions of $a$ with $b$ parts that contain $s$ staircases $1^+2^+3^+\cdots m^+$.

The main theorem arrived at by the end of the paper consists in establishing a formula for the generating function $F$ defined in equation (\ref{f36}). We state it here for completeness:
\begin{align*}
F&= \frac {k_{m}-\frac {qx^{m}y}{1-x}k_{m-1}}{(1-q)x^{\binom {m+1}{2}}\left(\frac{y}{1-x}\right)^{m}+\frac{1-x-xy}{1-x}\left(k_{m}-\frac{qx^{m}y}{1-x}k_{m-1}\right)},
\end{align*}
where $k_{m}=\sum_{\j=0}^{m-1}x^{mj-\binom {j}{2}}\left(\frac {y}{1-x}\right)^{j}$.
Prior to this main theorem, several lemmas present a set of recursions which are used in proving this result.

\section{Proofs}
\subsection{Warmup: compositions containing words of the form $1^+2^+$ or $1^+2^+3^+$}
Consider words which are of the form $1^+2^+$; i.e., words of two parts adjacent to each other from left to right with the first being a letter $>0$ and the second being a letter $>1$.

We let $F$ be the generating function for all words; $F_{a}$ be the generating function for all words starting with the letter $a$ and in general $F_{a_{1}a_{2}\cdots a_{n}}$ be the gf (generating function) for words starting with the letters $a_{1}a_{2}\cdots a_{n}$.
So by definition
\begin{equation}\label{f1}
F=1+\sum_{a\geq1}F_{a}.
\end{equation}
And we have the following recurrence:
\begin{equation}\label{f2}
F_{a}=x^ay+F_{a1}+F_{a2}+F_{a3}+\cdots
\end{equation}
Now $F_{a1}=x^ayF_{1}$ and $F_{ab}=qx^ayF_{b}$ for $b>1$. So $F_{a}=x^ay+F_{1}+qF_{2}+qF_{3}+\cdots$. Thus for all $a\ge1$, we have $F_{a}=x^ay(1-q)(1+F_{1})+qx^ayF$. As the second part of our warmup, we now examine the pattern $1^+2^+3^+$, i.e., we focus on compositions which contain this word sequence.

Extracting part of the first letter, we have
\begin{equation}\label{f5}
F_{a}=x^{a-1}F_{1}.
\end{equation}
From equation ({\ref{f1}}),
\begin{equation}\label{f6}
F=1+\sum_{a\geq1}F_{a}=1+\frac{1}{1-x}F_{1}.
\end{equation}
Also \begin{alignat}{1}
F_{1}&=xy+(F_{11}+F_{12}+F_{13}+\cdots)\nonumber\\
&=xy+xyF_{1}+F_{12}+xF_{12}+x^2F_{12}+\cdots)\nonumber\\
&=xy+xyF_{1}+\frac{1}{1-x}F_{12},\label{f7}
\end{alignat}
where
\begin{alignat}{1}
F_{12}&=x^3y^2+F_{121}+F_{122}+(F_{123}+\cdots)\nonumber\\
&=x^3y^2+x^3y^2F_{1}+x^2yF_{12}+(qx^3yF_{12}+qx^4yF_{12}+\cdots)\nonumber\\
&=x^3y^2+x^3y^2F_{1}+x^2yF_{12}+\frac{qx^3y}{1-x}F_{12}\label{f8}.
\end{alignat}

The last three equations have three unknowns $F, F_{1},$ and $F_{12}$ which we can solve for F using Cramer's rule. However, instead, we try the general pattern.

\subsection{The general pattern $1^+2^+3^+\cdots m^+$}
As before, $F_{a}=x^{a-1}F_{1}$ and
\begin{equation}\label{f7}
F=1+\sum_{a\geq1}F_{a}=1+\frac{1}{1-x}F_{1}.
\end{equation}
Now \begin{alignat}{1}
F_{1}&=xy+(F_{11}+F_{12}+F_{13}+\cdots)\nonumber\\
&=xy+xyF_{1}+F_{12}+xF_{12}+x^2F_{12}+\cdots)\nonumber\\
&=xy+xyF_{1}+\frac{1}{1-x}F_{12}\label{f10}
\end{alignat}
and
\begin{alignat}{1}
F_{12}&=x^3y^2+F_{121}+F_{122}+(F_{123}+\cdots)\nonumber\\
&=x^3y^2+x^3y^2F_{1}+x^2yF_{12}+(F_{123}+xF_{123}+x^2F_{123}+\cdots)\nonumber\\
&=x^3y^2+x^3y^2F_{1}+x^2yF_{12}+\frac{1}{1-x}F_{123}\label{f11}.
\end{alignat}
Next, by a similar process
\begin{alignat}{1}
F_{123}&=x^6y^3+x^6y^3F_{1}+x^5y^2F_{12}+x^3yF_{123}+\frac{1}{1-x}F_{1234}.\label{f12}
\end{alignat}
Proceeding in this way, we obtain in general for all $j\le{m-1}$
\begin{align}
F_{12\cdots j}&=x^{\binom{j+1}{2}}y^j+x^{\binom{j+1}{2}-\binom{1}{2}}y^{j}F_{1}+x^{\binom{j+1}{2}-\binom{2}{2}}y^{j-1}F_{12}\notag\\
&+x^{\binom{j+1}{2}-\binom{3}{2}}y^{j-2}F_{123}+\cdots+x^{\binom{j+1}{2}-\binom{j}{2}}yF_{12\cdots j}+\frac{1}{1-x}F_{12\cdots j+1}.\label{f13}
\end{align}
with
\begin{alignat}{1}
F_{12\cdots m}=qx^myF_{12\cdots m-1}\label{f14}.
\end{alignat}
To simplify the presentation we put $z=\frac{-1}{1-x}$. Now, we rewrite equations
\eqref{f8}-\eqref{f14} in matrix form. So we first define the matrix $\bf A$ as
{\tiny\begin{equation*}
\left(\begin{array}{llllllll}
{1}&{z}&{0}&{0}&\cdots&\cdots&\cdots&{0}\\
{0}&1-x^{\binom{2}{2}-\binom{1}{2}}y&{z}&{0}&\cdots&\cdots&\cdots&{0}\\
{0}&-x^{\binom{3}{2}-\binom{1}{2}}y^2&1-x^{\binom{3}{2}-\binom{2}{2}}y&{z}&\cdots&&&{0}\\
\vdots&&&&&&&\vdots\\
{0}&-x^{\binom{m-1}{2}-\binom{1}{2}}y^{m-2}&-x^{\binom{m-1}{2}-\binom{2}{2}}y^{m-3}&-x^{\binom{m}{2}-\binom{3}{2}}y^{m-4}&\cdots
&-x^{\binom{m-1}{2}-\binom{m-2}{2}}y&{z}&{0}\\
{0}&-x^{\binom{m}{2}-\binom{1}{2}}y^{m-1}&-x^{\binom{m}{2}-\binom{2}{2}}y^{m-2}&-x^{\binom{m}{2}-\binom{3}{2}}y^{m-3}&\cdots
&-x^{\binom{m}{2}-\binom{m-2}{2}}y^2&1-x^{\binom{m}{2}-\binom{m-1}{2}}y&{z}\\
{0}&{0}&{0}&{0}&\cdots&{0}&-qx^{m}y&{1}\\
\end{array}\right)
\end{equation*}}
and
${\bf C}$ to be the vector $\left(x^{\binom{1}{2}},x^{\binom{2}{2}}y,x^{\binom{3}{2}}y^2,\cdots,x^{\binom{m-1}{2}}y^{m-2},x^{\binom{m}{2}}y^{m-1},0\right)^T$.
Then the matrix form of our equations is $\bf A\bf X=\bf C$ where it is the first entry of matrix $\bf X$ (the matrix of variables from equations \eqref{f8}-\eqref{f14}) that is our required generating function $F$. So defining {\bf{B}} as the matrix obtained from the above matrix {\bf A} by replacing its first column with the entries from $\bf C$; i.e.
{\tiny\begin{equation*}
\left(\begin{array}{lllllll}
{x^{\binom{1}{2}}}&{z}&{0}&\cdots&\cdots&\cdots&{0}\\
{x^{\binom{2}{2}}y}&1-x^{\binom{2}{2}-\binom{1}{2}}y&{z}&\cdots&\cdots&\cdots&{0}\\
{x^{\binom{3}{2}}y^{2}}&-x^{\binom{3}{2}-\binom{1}{2}}y^2&1-x^{\binom{3}{2}-\binom{2}{2}}y&\cdots&&&{0}\\
\vdots&&&&&&\vdots\\
{x^{\binom{m-1}{2}}y^{m-2}}&-x^{\binom{m-1}{2}-\binom{1}{2}}y^{m-2}&-x^{\binom{m-1}{2}-\binom{2}{2}}y^{m-3}&\cdots
&-x^{\binom{m-1}{2}-\binom{m-2}{2}}y&{z}&{0}\\
{x^{\binom{m}{2}}y^{m-1}}&-x^{\binom{m}{2}-\binom{1}{2}}y^{m-1}&-x^{\binom{m}{2}-\binom{2}{2}}y^{m-2}&\cdots
&-x^{\binom{m}{2}-\binom{m-2}{2}}y^2&1-x^{\binom{m}{2}-\binom{m-1}{2}}y&{z}\\
{0}&{0}&{0}&\cdots&{0}&-qx^{m}y&{1}\\
\end{array}\right).
\end{equation*}}
By Cramer's rule, we obtain
\begin{equation}\label{f30}
F=\frac {\det{\bf{B}}}{\det{\bf{A}}}.
\end{equation}

\subsection{Equations for $\det\bf A$ and $\det\bf B$ in a form that can be solved recursively}
Define the $m$x$m$ matrix ${\bf {N_{m}}}$, to be the first $m$ rows and columns of the $(m+1)$x$(m+1)$ matrix {\bf {A}}, but where the first column of {\bf{A}} has initially been replaced by the first $m$ entries of {\bf C}.
To simplify the notation further, we let $w_{ij}=x^{\binom{i}{2}-\binom{j}{2}}y^{i-j}$ and so explicitly written out,
\begin{equation*}\label{f18}
\bf {N_{m}}:=\left(\begin{array}{llllll}
x^{\binom{1}{2}}y^0&z&0&0&\cdots&0\\
x^{\binom{2}{2}}y&1-w_{21}&z&0&&\vdots\\
x^{\binom{3}{2}}y^2&-w_{31}&1-w_{31}&z&&\\
\vdots&&&&&\vdots\\
x^{\binom{m}{2}}y^{m-1}&-w_{m1}&\cdots&&\cdots&1-w_{m1}\\
\end{array}\right).
\end{equation*}
By cofactor expansions (initially along the last row of {\bf B}), we obtain
\begin{equation}\label{f32}
\det{\bf B}=\det{\bf N_{m}}+zqx^{m}y\det{\bf N_{m-1}}.
\end{equation}

And let $\bf C_{m-1}$ be the $(m-1)$x$(m-1)$ matrix obtained by deleting the first row and column of $\bf {N_{m}}$. So, for example,
\begin{equation*}\label{f37}
\bf {C_{4}}=\left(\begin{array}{llll}
1-w_{21}&z&0&0\\
-w_{31}&1-w_{32}&z&0\\
-w_{41}&-w_{42}&1-w_{43}&z\\
-w_{51}&-w_{52}&-w_{53}&1-w_{54}\\
\end{array}\right).
\end{equation*}
By employing cofactor expansions (also, initially along the last row of ${\bf A}$), we see that
\begin{equation}\label{f31}
\det{\bf A}=\det{\bf C_{m-1}}+zqx^{m}y\det{\bf C_{m-2}}.
\end{equation}
Again, by employing co-factor expansions along the last row of $\bf {C_{4}}$, we see that
\begin{equation*}
\det{\bf C_{4}}=(1-w_{54})\det{\bf C_{3}}+zw_{53}\det{\bf C_{2}}-w_{52}z^2\det{\bf C_{1}}+w_{51}z^3\det{\bf C_{0}},\label{f20}
\end{equation*}
where $\det{\bf C_{0}}:=1$. In general, a cofactor expansion along the last row of ${\bf C_{m}}$ yields for $m\ge1$
\begin{equation*}
\det{\bf C_{m}}=(1-w_{m+1m})\det{\bf C_{m-1}}+\sum_{\j=1}^{m-1}(-1)^{m-1-j}w_{m+1j}z^{m-j}\det{\bf C_{j-1}}.
\end{equation*}
Once again making the replacement $w_{ij}=x^{\binom{i}{2}-\binom{j}{2}}y^{i-j}$, we have for $m\ge1$
\begin{equation}\label{f21}
\det{\bf C_{m}}=(1-x^{m}y)\det{\bf C_{m-1}}+\sum_{\j=1}^{m-1}(-1)^{m-1-j}x^{\binom{m+1}{2}-\binom{j}{2}}y^{m+1-j}z^{m-j}\det{\bf C_{j-1}}.
\end{equation}
Dropping $m$ by $1$ and multiplying this equation by $-x^{m}yz$, we obtain
\begin{align}\label{f22}
&-x^{m}yz\det{\bf C_{m-1}}\notag\\
&\quad=-x^{m}yz(1-x^{m-1}y)\det{\bf C_{m-2}}+\sum_{\j=1}^{m-2}(-1)^{m-1-j}x^{\binom{m+1}{2}-\binom{j}{2}}y^{m+1-j}z^{m-j}\det{\bf C_{j-1}}.
\end{align}
By subtracting \eqref{f22} from \eqref{f21}, we obtain
\begin{align*}
&\det{\bf C_{m}}+x^{m}yz\det{\bf C_{m-1}}\\
&\qquad=(1-x^{m}y)\det{\bf C_{m-1}}+x^{m}yz(1-x^{m-1}y)\det{\bf C_{m-2}}+x^{2m-1}y^2z\det{\bf C_{m-2}}.
\end{align*}
Simplifying,
\begin{equation}\label{f23}
\det{\bf C_{m}}=(1-x^{m}y(1+z))\det{\bf C_{m-1}}+x^{m}yz\det{\bf C_{m-2}},
\end{equation}
where $\det{\bf C_{-1}}:=1$;  $\det{\bf C_{0}}=1$;  $\det{\bf C_{1}}=1-xy=1-w_{21}$.

For ease of notation in the remainder of the paper, we abbreviate $\det{\bf C_{m}}$ as  $C_{m}$, and define the generating function $C(t)=\sum_{m\geq0}C_{m}t^{m}$.
By multiplying equation \eqref{f23} by $t^{m}$ and then summing from $1$ to infinity, we obtain
\begin{equation*}
C(t)-1=tC(t)-(1+z)xytC(xt)+x^2yt^2zC(xt)+xyzt.
\end{equation*}
Therefore
\begin{equation}
 C(t)=\frac {1+xyzt} {1-t}-xytC(xt)\frac{1+z(1-xt)}{1-t}.
\end{equation}
Again to simplify the notation, substitute $f(t):=\frac {1+xyzt} {1-t}$ and $\varphi(t):=-xyt\frac{1+z(1-xt)}{1-t}$,
and iterate the previous equation to obtain:
\begin{alignat}{1}
C(t)=f(t)+\varphi(t)C(xt)=f(t)+\varphi(t)f(xt)
+\varphi(t)\varphi(xt)C(x^2t).
\end{alignat}
Repeatedly iterating (assuming $|x|<1$), we obtain
\begin{align*}
C(t)&=\sum_{j\geq0}f(x^{j}t)\prod_{\i=0}^{j-1}\varphi({x^{i}t})\\
&=\sum_{j\geq0}(-1)^{j}\frac {1+x^{j+1}yzt}{1-x^{j}t}x^{\binom{j+1}{2}}y^{j}t^{j} \prod_{\i=0}^{j-1}\frac{1+z(1-x^{i+1}t)}{1-x^{i}t}.
\end{align*}
Recall that $z=\frac {-1}{1-x}$ which implies $1+z=\frac {-x}{1-x}$. Therefore,
\begin{align*}
C(t)&=\sum_{j\geq0}(-1)^j(1+x^{j+1}yzt)x^{\binom{j+1}{2}}y^{j}t^{j}\frac {\prod_{\i=1}^{j}(1-\frac {zx^{i}t}{1+z})}{\prod_{\i=0}^{j}(1-x^{i}t)}(1+z)^{j}\\
&=\sum_{j\geq0}(-1)^j(1+x^{j+1}yzt)x^{\binom{j+1}{2}}y^{j}t^{j}(\frac {-x}{1-x})^{j}\frac {\prod_{\i=0}^{j-1}(1-x^{i}t)}{\prod_{\i=0}^{j}(1-x^{i}t)}\\
&=\sum_{j\geq0}\frac {(1+x^{j+1}yzt)x^{\frac {j(j+3)}{2}}y^{j}t^{j}}{(1-x)^{j}(1-x^{j}t)}.
\end{align*}
For further notational simplification, we let
\begin{equation*}
f_{j}=\frac {(1+x^{j+1}yzt)x^{\frac {j(j+3)}{2}}y^{j}t^{j}}{(1-x)^{j}(1-x^{j}t)}.
\end{equation*}
Finally, substituting for the remaining $z$ as above and using partial fractions
\begin{align*}
f_{j}&=\frac {x^{1+\frac{j(j+3)}{2}}y^{j+1}t^{j}}{(1-x)^{j+1}}+\frac {x^{\frac{j(j+3)}{2}}y^{j}(1-x-xy)t^{j}}{(1-x)^{j+1}(1-x^{j}t)}\\
&=\frac {x^{1+\frac{j(j+3)}{2}}y^{j+1}t^{j}}{(1-x)^{j+1}}+\frac {x^{\frac{j(j+3)}{2}}y^{j}(1-x-xy)t^{j}}{(1-x)^{j+1}}\sum_{k\geq0}x^{jk}t^{k}.
\end{align*}
Hence the $m$th coefficient of $C(t)$ is given by
\begin{align*}
C_{m}&=\frac {x^{\binom{m+2}{2}y^{m+1}}}{(1-x)^{m+1}}+\sum_{j=0}^{m}\frac {x^{\frac{j^{2}+3j}{2}-j^{2}+jm}y^{j}(1-x-xy)}{(1-x)^{j+1}}
\end{align*}
So, we obtain the following lemma.
\begin{lemma}\label{lem1}
The determinants $C_{m}$ of the matrices obtained from $\bf{N_{m+1}}$ (see equation (\ref{f18})) by deleting its first row and column are given by
\begin{equation}
C_{m}=x^{\binom{m+2}{2}}\left(\frac {y} {1-x}\right)^{m+1}+\frac {1-x-xy} {1-x}\sum_{j=0}^{m} x^{(m+1)j-\binom{j}{2}}\left(\frac {y} {1-x}\right)^{j}.\label{f34}
\end{equation}
\end{lemma}

For initial cases, we have $\det{\bf {N_{1}}}=1$ and $\det{\bf {N_{2}}}=1-xy-zxy$.
By a cofactor expansion along the last row, we obtain for $m\ge2$
\begin{align}
\det{\bf {N_{m}}}&=(1-x^{m-1}y)\det{\bf {N_{m-1}}}\notag\\
&+\sum_{\j=1}^{m-2}(-1)^{m-j}x^{\binom {m}{2}-\binom{j}{2}}y^{m-j}z^{m-1-j}\det{\bf {N_{j}}}+(-1)^{m-1}x^{\binom {m}{2}}y^{m-1}z^{m-1}.\label{f19}
\end{align}
Dropping $m$ by $1$ and multiplying this equation by $-x^{m-1}yz$ (a similar process to that used in a previous section), we obtain for $m\ge3$
\begin{align}
-x^{m-1}yz\det{\bf N_{m-1}}&=-x^{m-1}yz(1-x^{m-2}y)\det{\bf N_{m-2}}\notag\\
&+\sum_{\j=1}^{m-3}(-1)^{m-j}x^{\binom{m}{2}-\binom{j}{2}}y^{m-j}z^{m-1-j}\det{\bf N_{j}}
+(-1)^{m-1}x^{\binom {m}{2}}y^{m-1}z^{m-1}.\label{f25}
\end{align}
Subtracting \eqref{f25} from \eqref{f19}, we obtain
\begin{align*}
&\det{\bf N_{m}}+x^{m-1}yz\det{\bf N_{m-1}}\\
&=(1-x^{m-1}y)\det{\bf N_{m-1}}+x^{m-1}yz(1-x^{m-2}y)\det{\bf N_{m-2}}+x^{2m-3}y^2z\det{\bf N_{m-2}}\\
&=(1-x^{m-1}y)\det{\bf N_{m-1}}+x^{m-1}yz\det{\bf N_{m-2}}.
\end{align*}
Hence for $m\ge2$,
\begin{equation}\label{f26}
\det{\bf N_{m}}=(1-x^{m-1}y(1+z))\det{\bf N_{m-1}}+x^{m-1}yz\det{\bf N_{m-2}}
\end{equation}
with $\det{\bf N_{0}}=0$ and $\det{\bf N_{1}}=1$.

For the rest of the paper we simplify matters by abbreviating $N_m:=\det{\bf N_{m}}$ and now define the generating function $N(t)=\sum_{m\geq0}N_mt^m$. By multiplying equation \eqref{f26} by $t^{m}$, summing from $1$ to infinity, we obtain
\begin{equation*}
N(t)-t=tN(t)-y(1+z)tN(xt)+xyzt^{2}N(xt)
\end{equation*}
with $N_{-1}:=0$.
Hence
\begin{equation}\label{f28}
 N(t)=\frac {t} {1-t}+\frac {xyzt^{2}-y(1+z)t}{1-t}N(xt).
\end{equation}
Repeatedly iterating \eqref{f28} on $t$ (while recalling that $z=\frac {-1} {1-x}$, and assuming $|x|<1$), we obtain
\begin{align*}
N(t)&=\sum_{j\geq0}\frac {x^{j}t}{1-x^{j}t} \prod_{\i=0}^{j-1}\frac {yx^{i}t(\frac{-x^{i+1}t}{1-x}+\frac {x} {1-x})}{1-x^{i}t}\\
&=\sum_{j\geq0}\frac {x^{j}t}{1-x^{j}t} \prod_{\i=0}^{j-1}\frac {yx^{i}t}{1-x}\\
&=\sum_{j\geq0}\frac {x^\frac{j^{2}+3j}{2}y^{j}t^{j+1}}{(1-x^{j}t)(1-x)^{j}} .
\end{align*}
Thus, we have our final lemma.
\begin{lemma}\label{lem2}
With $N_{m}:=\det{\bf{N_{m}}}$ (see \eqref{f18})
\begin{equation}\label{f34}
N_{m}=[t^{m}]N(t)=\sum_{j=0}^{m-1}x^{mj-\binom {j}{2}}\left(\frac {y}{1-x}\right)^{j}.
\end{equation}
\end{lemma}

\subsection{The generating function $F$}
Finally, apply \eqref{f32} and \eqref{f31} to \eqref{f30}. Then, use lemma \ref{lem1} and lemma \ref{lem2}, to obtain:

\begin{theorem}\label{th1}
The generating function $F=\sum_{n\geq1; b\geq1; s\geq 0}n(a,b,c)x^{a}y^{b}q^{s}
$ for the number of staircases $1^+2^+3^+\cdots m^+$ (tracked by the exponent of variable $q$) contained in particular compositions (of $a$ with $b$ parts) is given by
\begin{equation}
F= \frac {N_{m}-\frac {qx^{m}y}{1-x}N_{m-1}}{(1-q)x^{\binom {m+1}{2}}\left(\frac{y}{1-x}\right)^{m}+\frac{1-x-xy}{1-x}\left(N_{m}-\frac{qx^{m}y}{1-x}N_{m-1}\right)}\label{f35}.
\end{equation}
\end{theorem}

For example, Theorem \ref{th1} with $q=1$ yields $F_{q=1}=\frac{1-x}{1-x-y}$, which is the generating function for the number of compositions of $n$ with exactly $m$ parts (see \cite{Mb1}).

By differentiating the generating function $F$ with respect to $q$ and then substituting $q=1$, we obtain
\begin{align*}
\frac{dF}{dq}\mid_{q=1}&=
\frac{x^{\binom {m+1}{2}}\left(\frac{y}{1-x}\right)^{m}}{\frac{(1-x-xy)^2}{(1-x)^2}\left(\sum_{j=0}^{m-1}x^{mj-\binom {j}{2}}\left(\frac {y}{1-x}\right)^{j}
-\sum_{j=1}^{m-1}x^{mj-\binom {j}{2}}\left(\frac {y}{1-x}\right)^{j}\right)}\\
&=\frac{x^{\binom {m+1}{2}}y^m}{(1-x-xy)^2(1-x)^{m-2}}\\
&=\frac{x^{\binom {m+1}{2}}}{(1-x)^m}\sum_{j\geq0}(j+1)\frac{x^jy^{m+j}}{(1-x)^j}
\end{align*}
Next, we extract coefficients; firstly of $[y^{l}]$  to obtain
$$(\ell-m+1)\frac{x^{\ell+\binom {m}{2}}}{(1-x)^\ell}=(\ell-m+1)\sum_{j\geq0}\binom{\ell+j-1}{j}x^{\ell+j+\binom {m}{2}},$$

and then of $[x^n]$
which leads to the following result.
\begin{corollary}
The total number of staircases $1^+2^+3^+\cdots m^+$ in all compositions of $n$ with exactly $\ell$ parts is given by
$$(\ell-m+1)\binom{n-1-\binom {m}{2}}{\ell-1}.$$
\end{corollary}


%
%
%
%
%
%
%
%
%
%
%

\end{document}